\documentclass[11pt]{amsart}
\usepackage{amsmath,amssymb,amsthm}
\usepackage[margin=1.1in]{geometry}
\usepackage{tikz,enumitem}
\usepackage{hyperref}
\usepackage{esvect}
\usepackage{float}
\usepackage{enumitem}
\usepackage{mathtools}
\usepackage{xcolor}

\newcommand{\Q}{\ensuremath{\mathbb Q}}

\newtheorem{theorem}{Theorem}
\newtheorem*{reptheorem}{Theorem \ref{mainTheorem} (Alternate version)}

\newtheorem{example}[theorem]{Example}
\newtheorem{corollary}[theorem]{Corollary}
\newtheorem{lemma}[theorem]{Lemma}
\theoremstyle{definition}
\newtheorem{definition}[theorem]{Definition}
\newtheorem{remark}[theorem]{Remark}

\title{The MZ-property for three-term differences}

\author{
  Hajrudin Fejzi\'c and
  Chris Freiling \\\\ California State University, San Bernardino
}
\date{} 

\begin{document}

\begin{abstract} We study the Marcinkiewicz--Zygmund property for finite differences. A difference
\[
D(f,h)=\sum_{k=1}^n a_k f(b_kh)
\]
is said to have the MZ-property if, for every measurable function \(f\),
\[
D(f,h)=o(h) \quad\text{and}\quad f(h)=o(1)
\]
imply
\[
f(h)=o(h).
\]
In this paper, we characterize all three-node differences
\[
D(f,h)=af(qh)+bf(rh)+cf(sh)
\]
having the MZ-property. In particular, for irrational differences with
distinct nodes, we obtain explicit necessary and sufficient conditions
in terms of the coefficients \(a,b,c\) and the nodes \(q,r,s\).
As an application, our results resolve an open problem concerning a
second-order Riemann difference.
\end{abstract}
\footnotetext{\textbf{2020 Mathematics Subject Classification.}
Primary 26A24; Secondary 26A27.
\noindent\textbf{Keywords.}
differences, Riemann derivative, MZ-property.}

\maketitle

We consider three-term functional differences of the form
\begin{equation}\label{Diff3}
D(f,h) = af(qh)+bf(rh)+cf(sh),
\end{equation}
where $h$ is a real variable and $f$ is an arbitrary complex-valued function defined for real 
numbers in a neighborhood of 0. 
Each such difference $D$ has fixed parameters.  These include nonzero coefficients, $a,b,c$, which may be real or complex, and distinct nodes, $q,r,s$, which are nonzero real numbers.
We consider three possible limit properties of $f$ as $h \to 0$:
\begin{itemize}
\renewcommand\labelitemi{} 
    \item A. $D(f,h) = o(h)$
    \item B. $f(h) = o(1)$
    \item C. $f(h) = o(h)$
\end{itemize}
Property C trivially implies both A and B, for any choice of parameters.  If A and B together imply C then the difference $D(f,h)$ is said to have the \emph{MZ-property}.

A case requiring special consideration is the following.
\begin{definition}
    A difference~\eqref{Diff3} is called \emph{rational} if there are three integers $i,j,k$, not all zero, such that $i\log|q|+ j\log|r| + k\log|s|=i+j+k=0$.  Otherwise, it is called \emph{irrational}.
\end{definition}

A characterization of rational differences with the MZ-property is contained in the results of ~\cite{HF} and will be discussed below.
In this paper we characterize the irrational three-term differences that have the MZ-property.  Our main result is the following:
\begin{theorem}\label{mainTheorem}
    Let $D(f,h) = af(qh)+bf(rh)+cf(sh)$ be an irrational difference with distinct nodes,
    $q,r,s$, where 
    $|q|\le|r|\le|s|$.
    Then $D$ has the MZ-property if and only if at least one
    of the following conditions holds
    \begin{itemize}
\item $|a|+|b|\leq |c|$
\item $|br| + |cs| < |aq|$
\item $|a|+|c|\leq |b|$ and $|cs| + |aq|< |br|$
    \end{itemize}
\end{theorem}

Since the MZ-property is invariant under scaling either the coefficients $\{a,b,c\}$ or the nodes $\{q,r,s\}$,
we will often rewrite the difference in \emph{standard form}, with $q=c=1$, $1\le |r| \le |s|$, and
\begin{equation}\label{Standard}
D(f,h) = af(h)+bf(rh)+f(sh).
\end{equation}
Then the definition of \emph{rational} takes on the more familiar form:
\begin{definition}
    A difference~\eqref{Standard} is called \emph{rational} if and only if $\frac{\log|r|}{\log|s|}\in \Q$.  
\end{definition}
Theorem~\ref{mainTheorem} can then be re-written in a less symmetrical, but equivalent form:
\begin{reptheorem}
    Let $D(f,h) = af(h)+bf(rh)+f(sh)$ be an irrational difference with 
    distinct nodes, $1,r,s$ where $1\le |r|\le |s|$.  Then $D$ has the MZ-property if and only if at least one
    of the following conditions holds
    \begin{itemize}
\item $|a|+|b|\leq 1$
\item $|br| + |s| < |a|$
\item $|a|+1\leq |b|$ and $|s| + |a|< |br|$.
    \end{itemize}
\end{reptheorem}

\section{Real differences}

The original motivation for the MZ-property was to study smoothness properties of real functions.
However, as shown in~\cite{HF}, this property cannot be properly studied without also including
complex numbers.  Hence, we allow complex numbers for the coefficients, $a,b,c$, as well as
the output of the arbitrary function $f$.  In this section we will justify 
the use of complex numbers for the case that we are only interested in real numbers.

\begin{definition}
    A difference~\eqref{Standard} is called \emph{real} if the coefficients, $a,b,c$ are all real numbers.  
\end{definition}

\begin{definition}
    A difference~\eqref{Standard} has the \emph{real} MZ-property if A and B imply C for all real-valued functions.
\end{definition}
\begin{theorem}\label{realMZ}
A real difference~\eqref{Standard} has the MZ-property if and only if it has the real MZ-property.
\end{theorem}
\begin{proof}
First, suppose that a real difference, $D$, has the MZ-property.  Then properties A and B together imply C for any function $f:\mathbb R \to \mathbb C$.  Trivially,  A and B together imply C for any function $f:\mathbb R \to \mathbb R$.  

Next, suppose that a real difference fails the MZ-property.  Then there is a function $f:\mathbb R \to \mathbb C$ such that 
$D(f,h) = o(h)$,
$f(h) = o(1)$,
but $f(h)$ is not $o(h)$.  Let $r(h) = Re(f)$, and $c(h) = Im(f)$,
so that $r:\mathbb R \to \mathbb R$ and $c:\mathbb R \to \mathbb R$.
Since $f=o(1)$, we have both $r=o(1)$ and $c=o(1)$.  
Since $D(f,h) = o(h)$, we have $Re(D(f,h)) = o(h)$ and $Im(D(f,h)) = o(h)$.
But since $D$ is a real difference, $Re(D(f,h)) = D(Re(f),h) = D(r,h)$,
and $Im(D(f,h)) = D(Im(f),h) = D(c,h)$.
Finally, if $r(h)$ and $c(h)$ were both $o(h)$, then we would have $f(h)=o(h)$.
By choice of $f$, this does not happen.  So, at least one of $r(h), c(h)$
witnesses that $D$ does not have the real MZ-property.
\end{proof}

A consequence of this is that Theorem~\ref{mainTheorem} also characterizes the real MZ-property
for real differences. 
For the remainder of this paper, we will only discuss the complex variant.

\section{Background}
When an ordinary derivative is not available, one may consider generalized derivatives such as the symmetric derivative
\[
f_s'(c)=\lim_{h\to0}\frac{f(c+h)-f(c-h)}{2h}.
\]
The symmetric derivative is a generalization of ordinary differentiation since if $f'(c)$ exists then so does 
$f_s'(c)$ and they are equal.  Furthermore, it is a proper
generalization since 
\[
f(h)=|h|
\]
is symmetrically differentiable at \(0\), but not ordinarily differentiable there.

More generally, derivatives of the form 
\begin{equation}\label{GenDef}
   D'f(c)=\lim_{h\to0}\frac{\sum_{k=1}^n a_kf(c+b_kh)}{h}. 
\end{equation}
can be considered, where the coefficients $a_i$ are real numbers.

To generalize the ordinary derivative, we need the derivative of a constant to be 0
and the derivative of the identity function to be 1.  
That is, we need
\begin{equation}\label{GenReq}
\sum_{k=1}^n a_k=0
\qquad\text{and}\qquad
\sum_{k=1}^n a_kb_k=1.
\end{equation}
It is not hard to see that these two conditions guarantee that $D'f(c)$ is a generalization of the ordinary derivative.  That is, under these two conditions, whenever $f'(c)$  exists, so does $D'f(c)$, and $D'f(c)=f'(c)$. For this reason, derivatives that satisfy 
these two conditions are called \emph{Riemann derivatives}.

A Riemann derivative, however, might not yield anything new; it may be equivalent to ordinary 
differentiation.  For example, let
\[
S'f(c)=\lim_{h\to0}\frac{2f(c+h)+f(c-h)-3f(c)}{h}.
\]
If \(S'f(c)\) exists, then $f'(c)$ exists, and $f'(c) = S'f(c)$.
So, unlike the symmetric derivative $S'$ is not a proper generalization of the ordinary derivative.

It is natural to wonder, therefore, which Riemann derivatives are
equivalent to the ordinary derivative and which are proper generalizations.
Ash, Catoiu, and Cs\"ornyei \cite{ACC} gave a complete answer to this question.
Specifically, a derivative 
\eqref{GenDef} is equivalent to the ordinary derivative if and
only if it can be written in the form
\[
R'f(c)=
\lim_{h\to 0}
\frac{
af(c+qh)+af(c-qh)-2af(c) +bf(c+rh)-bf(c-rh)
}{2h},
\]
where $aq\neq0$ and $br = 1$. The derivative $S'$ above is obtained from $R'$ by setting $a=3$ and $q=r=b=1.$

This raises the question: What if we add continuity as an additional assumption?
That is, when does

\begin{equation}
D'f(c)=k
\qquad\text{and}\qquad
\lim_{h\to 0}f(c+h)=f(c)
\end{equation}
imply
\[
f'(c)=k?
\]
To study this question, we
assume, without loss of generality, that $c=0$ and subtract an affine function so that $f(c)=k=0$. This leads us to the following definition,
named after  Marcinkiewicz and Zygmund  since they are believed to be the first to consider this property (see \cite{MZ}).

\begin{definition}
A finite difference 
\[
D(f,h)=\sum_{k=1}^n a_kf(b_kh),
\]
is said to have the \emph{MZ-property} if and only if for any real function $f$ defined in a neighborhood of 0, 
\begin{equation}\label{e1}
D(f,h)=o(h) \text{ and } f(h)=o(1)
\implies f(h)=o(h).\end{equation} 
\end{definition}

Note that although the MZ-property was motivated by the study of Riemann derivatives, it may have intrinsic and independent interest, so we make no tacit assumptions on the properties~\eqref{GenReq} in this paper.
Also note, as discussed above, that we have changed this definition to include complex coefficients and functions $f:\mathbb R \to \mathbb C$. Finally for simplicity, we have ruled out the possibility that $b_k=0$.  This is without loss of generality, since 
otherwise, $D(f,h)=o(h)$ and $f(h)=o(1)$ imply that $f(b_kh)=f(0)=0$.

\begin{remark}\label{r1}
In their study of $m$-th Riemann derivatives \footnote{An \(m\)-th
Riemann derivative is determined by a difference
$D(f,h)=\sum_{k=1}^{m+1}a_kf(b_kh),$
whose coefficients and nodes satisfy
$\sum_{k=1}^{m+1}a_kb_k^j=0,
\qquad j=0,1,\ldots,m-1,
$ and
$\sum_{k=1}^{m+1}a_kb_k^m=m!. $  The corresponding \(m\)-th Riemann derivative at $c$ is
\[
D^m f(c)=\lim_{h\to0}\frac{D(f_c,h)}{h^m},
\qquad f_c(x)=f(c+x).
\]
}, Marcinkiewicz and Zygmund needed, in place of \eqref{e1}, the implication
\[D(f,h)=o(h^m)
\quad\text{and}\quad
f(h)=o(h^{m-1})
\quad\Longrightarrow\quad
f(h)=o(h^m).
\]
Here \(m\) is a positive integer, and we call this the MZ-property of
order \(m\). However, the MZ-property of order \(m\) reduces to the MZ-property
defined above. Indeed, define
\[
D_1(g,h)=\sum_{k=1}^n a_kb_k^{m-1}g(b_kh).
\]
Then \(D\) has the MZ-property of order \(m\) if and only if \(D_1\)
has the MZ-property. To see this, suppose first that \(D\) fails the
MZ-property of order \(m\) for some \(f\), and set
\[
g(h)=\frac{f(h)}{h^{m-1}}.
\]
Then \(g(h)=o(1)\), while
\[
D(f,h)=h^{m-1}D_1(g,h).
\]
Thus \(D_1(g,h)=o(h)\), whereas \(g(h)\ne o(h)\), so \(D_1\) fails the
MZ-property. Conversely, if \(D_1\) fails the MZ-property for some
\(g\), then setting
\[
f(h)=h^{m-1}g(h)
\]
shows that \(D\) fails the MZ-property of order \(m\).
\end{remark}

The MZ-property has been studied for higher order derivatives and for various number of terms, (See \cite{ACF1}, \cite{ACF}, and \cite{CF1}.) though
a complete characterization has not yet appeared. The primary result of this study solves the case of 
\emph{geometric nodes}, defined as follows.

\begin{definition}
  Let $\lambda$ be a real number with $|\lambda|>1$. Then the finite difference  
  $$D(f,h)=\sum_{k=0}^n \left( a_kf(\lambda^k h) + c_kf(-\lambda^k h)\right)$$
  is said to have \emph{geometric nodes} with ratio $\lambda$. It has two \emph{characteristic polynomials},
  \begin{equation}
  P_+(z) = \sum_{k=0}^n (a_k+c_k) z^{k}
  \qquad\text{and}\qquad
  P_-(z) = \sum_{k=0}^n (a_k-c_k) z^{k}
  \end{equation}
\end{definition}
\noindent As an example, suppose $D(f,h)=5f(-2h)+3f(-h)-5f(0)-7f(2h)+f(2^3h)$.  Then $\lambda = 2$,the characteristic polynomials are $P_+(z) = z^3-2z+3$, and  $P_-(z) =z^3-12z-3$.

In  the case of geometric nodes, 
the MZ-property has been completely characterized in terms of the roots of the characteristic polynomials. 
\begin{theorem}(Fejzi\'c \cite{HF})\label{rationalTh}
    A difference with geometric nodes has the MZ-property if and only if neither of 
    its characteristic polynomials has a complex root $z$ in the annulus $1<|z|\leq|\lambda|$.
\end{theorem}
\noindent In the course of that investigation, it was necessary to allow for complex coefficients, explaining why
we retain this convention, even if our main interest is real-valued functions.

Any two-term difference can be written with geometric nodes after scaling.  Therefore, Theorem~\ref{rationalTh} completely characterizes the two-term case:
\begin{corollary}(Fejzi\'c \cite{HF})
Let $D(f,h) = af(h)+ f(qh)$ be a two-term difference with $1\le|q|$.  \newline
Case 1:
If $1=|q|$, i.e. $q=-1$ then $D(f,h)$ has the MZ-property if and only if 
$a \ne \pm 1$.\newline
Case 2:
If $1<|q|$ then $D(f,h)$ has the MZ-property if and only if
$a$ is not in the annulus $1<|a|\leq |q|$. 
\end{corollary}
\begin{proof}
Case 1: This can be deduced from Theorem~\ref{rationalTh}, but is also seen directly. Suppose $D(f,h)=af(h)+f(-h)=o(h)$. Then  
$aD(f,h)-D(f,-h)=(a^2-1)f(h)=o(h)$. If $a^2\ne 1$ then $f(h)=o(h)$. For counterexamples, if $a=1$
use $f(h)=h$ and if $a=-1$ use $f(h)=|h|$.

Case 2: Let $\lambda = |q|$.  Then, $D(f,h) = a f(\lambda^0 h) + f(\lambda^1 h)$ or $D(f,h) = a f(\lambda^0 h) + f(-\lambda^1 h)$.
In any event, the two characteristic polynomials are $a+z$ and $a-z$, with roots $z=\pm a$.
By Theorem~\ref{rationalTh}, the MZ-property fails if and only if $1<|a|\leq |q|$.
\end{proof}

For rational three-term differences in standard form~\eqref{Standard}, if
\[
\frac{\log|r|}{\log|s|}= \frac{k}{j},
\]
for integers $k\ge 0$ and $j\ge 1$ then letting $\lambda=|s|^{1/j}$ the difference can be written with geometric nodes
as $D(f,h) = af(h)+bf(\operatorname{sgn}(r)\lambda^k h)+f(\operatorname{sgn}(s)\lambda^j h)$, where $|\lambda|>1$.
Therefore, Theorem~\ref{rationalTh} applies, giving the following corollary.

\begin{corollary}(Fejzi\'c \cite{HF})\label{rat3}
Let $D(f,h) = af(h)+bf(rh)+f(sh)$ be a rational difference with $1 \le  |r| \le |s|$, and let $j>0,k\ge 0$ be integers such that
$\frac{\log|r|}{\log|s|}= \frac{k}{j}$.  Then $D$ has the MZ-property if and
only if the trinomials $z^j + bz^k + a$ and
$\operatorname{sgn}(s)z^j + \operatorname{sgn}(r)bz^k + a$
have no roots in the annulus $1<|z|\leq |s|^{1/j}$.
\end{corollary}

The next examples illustrate important differences between the rational and irrational cases.
Examples~\ref{Ex1} and \ref{Ex2} illustrate that the
MZ-property in the rational case can depend on the plus and minus signs attached to coefficients.  Examples~\ref{Ex1} and \ref{Ex3}  illustrate that it can also depend on plus and minus signs attached to the nodes.  This is in contrast to the irrational case, Theorem~\ref{mainTheorem},
where the MZ-property only depends on absolute values.
\begin{example}\label{Ex1}
The difference $D(f,h)=f(8h)+4f(2h)+3f(h)$ has only one distinct trinomial $x^3+4x+3$ whose roots are outside the annulus $1<|z|\leq 2$, so it has the MZ-property.
\end{example}
\begin{example}\label{Ex2}
The difference $D(f,h)= f(8h)-4f(2h)+3f(h)$ has only one distinct trinomial $x^3-4x+3$. 
It contains a root, $z$, where $1<|z|\leq 2$, so it does not have the MZ-property.
Indeed, here is a specific counterexample.  Let $f(h)=z^n$ if $h=2^n$ and $n\in \mathbb{Z}$.  Let $f(h)=0$ elsewhere.
\end{example}
\begin{example}\label{Ex3}
The difference $D(f,h)=f(8h)+4f(-2h)+3f(h)$ has two distinct trinomials, $x^3+4x+3$ and $x^3-4x+3$ and since $z=\frac{\sqrt{13}-1}{2}$ is a root of the second trinomial with $1<|z|\leq 2$, it fails MZ-property. Indeed, one can verify that $f(h)=\pm z^n$ for $h=\pm 2^n$ $n\in \mathbb{Z}$ and 0 elsewhere is a counterexample.
\end{example}

Accordingly, the main focus of the present paper is the case of irrational three-term differences.
Together, Theorem~\ref{mainTheorem} and Corollary~\ref{rat3} completely characterize the MZ-property for three-term differences.  But there is a big difference in these results.
Theorem~\ref{mainTheorem} solves the irrational case with a direct computation on the coefficients and nodes.
Corollary~\ref{rat3} solves the rational case, using the location of the complex roots of polynomial equations.  We do not know of a direct characterization for the rational case.
  In particular, Theorem~\ref{mainTheorem} does not hold for the rational case. 
  This is seen in Example~\ref{Ex1} where we have $a=3, b=4, r=2, s=8$. None of the three conditions of Theorem~\ref{mainTheorem} hold, and yet $D(f,h)$ has the MZ-property. Nevertheless, one direction of Theorem~\ref{mainTheorem} does hold in the rational case.
  
  \begin{theorem}\label{rationalTheorem}
    Let $D(f,h) = af(h)+bf(rh)+f(sh)$ be a rational difference with 
    distinct nodes, $1,r,s$ where $1\le |r|\le |s|$.  Then $D$ has the MZ-property if at least one of the following conditions holds
    \begin{itemize}
\item $|a|+|b|\leq 1$
\item $|br| + |s| < |a|$
\item $|a|+1\leq |b|$ and $|s| + |a|< |br|$.
    \end{itemize}
  \end{theorem}

We end this section with an application of Theorem~\ref{mainTheorem}
to an open problem posed in \cite{PMZ}. Consider the second-order Riemann difference
\[ D(f,h)=f(5h)-2f(2h)+f(-h).\]
We show that
\[
D(f,h)=o(h^2)
\quad\text{and}\quad
f(h)=o(h)
\]
imply
\[
f(h)=o(h^2).
\]
By Remark~\ref{r1}, it suffices to show that the difference
\[
D_1(g,h)=5g(5h)-4g(2h)-g(-h)
\]
has the MZ-property. Indeed, using the notation of
Theorem~\ref{mainTheorem}, the coefficients \(a=-1\), \(b=-4\),
and \(c=5\) satisfy the first condition of the theorem,
\[
|a|+|b|\leq |c|.
\]

The main goal now, is to establish Theorems~\ref{mainTheorem} and \ref{rationalTheorem}.
The rest of this paper is organized as follows. 
In Section~\ref{EdgeCases} we establish Theorems~\ref{mainTheorem} and
\ref{rationalTheorem} for the edge cases $r=-1$ and $r=-s$.  This is a mere convenience, allowing us to assume strict inequalities
$1<|r|<|s|$ for the remainder of the paper.
In Sections~\ref{sec2}, \ref{sec3}, and ~\ref{sec4}
we prove the MZ-property for the three respective conditions stated in Theorem~\ref{mainTheorem}
and Theorem~\ref{rationalTheorem}.
In Section~\ref{sec5} we show that if none of the three conditions hold, and $D(f,h)$ is irrational, then the MZ-property fails.

\section{The edge cases \texorpdfstring{$r= -1$ and $r=-s$}{$r= -1$ and $r=-s$}}\label{EdgeCases}

In this section we consider the difference
$$D(f,h) = af(h)+bf(rh)+f(sh)$$ where $1,r,s$ are distinct,
$1\le |r|\le |s|$, and one of the inequalities in not strict.  That is, either
$r= -1$ or $r=-s$.

The cases $r= -1$ and $r=-s$ do not occur in the irrational case.
Indeed, if $r=-1$ then 
$\frac{\log |r|}{\log |s|} = 0$ and if $r=-s$ then 
$\frac{\log |r|}{\log |s|} = 1$.  In both cases, $D(f,h)$ is rational.
The goal of this section, then, is to establish Theorem~\ref{rationalTheorem} 
for these two cases.

Let us first consider $r=-1$, so that
\begin{equation}\label{r=-1}D(f,h) = af(h)+bf(-h)+f(sh).
\end{equation}
with $|s|>1$.
\begin{corollary}\label{thr=-1}
A difference~\eqref{r=-1} has the MZ-property if and only if
neither $z=|a+b|$ nor $z=|a-b|$ is in the interval $1<|z|\le |s|$.
\end{corollary}
\begin{proof}
Apply Corollary~\ref{rat3}, with $k=0$ and $j=1$.
The two characteristic polynomials are
$z + (b + a)$ and
$\operatorname{sgn}(s)z + (a-b)$.  The absolute values of the respective roots 
are $|z|=|a+b|$ and $|z|=|a-b|$.
The result immediately follows from Corollary~\ref{rat3}.

\end{proof}
The following corollary is Theorem~\ref{rationalTheorem} when $r=-1$.
\begin{corollary}
If any of the following conditions hold
 \begin{itemize}
\item $|a|+|b|\leq 1$
\item $|b| + |s| < |a|$
\item $|a|+1\leq |b|$ and $|s| + |a|< |b|$,
    \end{itemize}
    then the difference~\eqref{r=-1} has the MZ-property.
\end{corollary}
\begin{proof}
Note that $$\max \{|a|-|b|, |b|-|a|\} \le \max \{|a-b|, |a+b|\} \le |a|+|b|.$$
So, if $|a|+|b|\le 1$ then  $|a+b|\le 1$ and $|a-b|\le 1$.  
If $|b| + |s| < |a|$, then $|s| < |a|-|b|$, so  $|s| < |a+b|$ and $|s| < |a-b|$. 
Finally, if $|s| + |a|< |b|$ then $|s| < |b|-|a|$, so  $|s| < |a+b|$ and $|s| < |a-b|$.
In all three cases, by Corollary~\ref{thr=-1}, $D(f,h)$ has the MZ-property.
\end{proof}

Next, consider $r=-s$, so that
\begin{equation}\label{r=-s}D(f,h) = af(h)+bf(-sh)+f(sh).
\end{equation}
with $|s|>1$.
\begin{corollary}\label{thr=-s}
A difference~\eqref{r=-s} has the MZ-property if and only if
$|a|$ is not in either of the intervals:
\begin{itemize}
\item $(|1+b|, |s||1+b|]$
\item $(|1-b|, |s||1-b|]$
    \end{itemize}
\end{corollary}
\begin{proof}
Apply Corollary~\ref{rat3}, with $k=1$ and $j=1$.
The two characteristic polynomials are
$$P_+(z)=(1+b)z + a \qquad \text{and} \qquad
P_-(z)=(\operatorname{sgn}(s)+ \operatorname{sgn}(-s)b)z+ a.$$  
Since the parameters are nonzero, if $b=-1$ then $P_+$ has no root.  Otherwise, it has a single root
with absolute value 
$|z|=\frac{|a|}{|1+b|}$.  In other words,  $P_+$ has no root, $z$, in the annulus $1<|z|\le |s|$ if and only if
$|a|$ is not in the interval $(|1+b|, |s||1+b|]$.
Similarly, $P_-$ has no root with $1<|z|\le |s|$ if and only if
$|a|$ is not in the interval $(|1-b|, |s||1-b|]$.
The result immediately follows by Corollary~\ref{rat3}.
\end{proof}
The following corollary is Theorem~\ref{rationalTheorem} when $r=-s$.
\begin{corollary}
If any of the following conditions hold
 \begin{itemize}
\item $|a|+|b|\leq 1$
\item $|bs| + |s| < |a|$
\item $|a|+1\leq |b|$ and $|s| + |a|< |bs|$,
    \end{itemize}
    then the difference~\eqref{r=-s} has the MZ-property.
\end{corollary}
\begin{proof}
If $|a|+|b|\leq 1$ then $|a|\le 1-|b|\le \min\{ |1+b|,|1-b|\}$.
If $|bs| + |s| < |a|$ then $|s|\max\{  |1+b|,|1-b|\}<|a|$.
If $|a|+1\le |b|$ then $|a|\le |b|-1 \le \min\{  |1+b|,|1-b|\} $.
In all three cases, the MZ-property follows by Corollary~\ref{thr=-s}.
\end{proof}
It remains to show Theorems~\ref{mainTheorem} and \ref{rationalTheorem} for
the case $1<|r|<|s|$.

\section{The case \texorpdfstring{$|a|+|b|\le 1$}{|a|+|b| <= 1}}\label{sec2}

In this section we consider
$$D(f,h) = af(h)+bf(rh)+f(sh)$$
where $1<|r|<|s|$  and $|a|+|b|\le 1$.

In the rational case, suppose $\frac{\log|r|}{\log|s|}= \frac{k}{j}$ for some
positive integers, \(j, k\). Then, this result follows from Corollary~\ref{rat3}.
In fact, no polynomial $P(z) = \pm z^j \pm bz^k \pm a$ has a root with $|z|>1$.
To see this, let $|z|>1$.  Then, 
\[
|P(z)| \ge |z^j| -|bz^k| - |a| > |z|^j(1-|b|) - |a| > 1-|b|-|a|\ge 0.
\]
Therefore, $P(z)$ has no zeros with $|z|>1$.
By Corollary~\ref{rat3},
$D(f,h)$ has the MZ-property. 

The proof below works directly
and does not require any rational relation between \(\log |r|\) and
\(\log |s|\).

\begin{theorem}\label{MZ}
Let \(1<|r|<|s|\), and let \(a,b\in\mathbb C\) satisfy
\[
|a|+|b|\le 1.
\]
Then
\[
Df(h)=af(h)+bf(rh)+f(sh)
\]
has the MZ-property
\end{theorem}

\begin{proof}
Assume $Df(h)=o(h)$ and $f(h)=o(1)$ and fix \(\varepsilon>0\). 
Then,
there exists \(\delta>0\) such that for all $t$ with \(0 < |t|<\delta\),
\[
|af(t)+bf(rt)+f(st)|\le \varepsilon |t| < \varepsilon |st|.
\]
Then,
\begin{equation}\label{stIneq}
    |f(st)| < |a|\left| f\left(t\right)\right| + |b|\left| f\left(rt\right)\right| + \varepsilon |st|.
\end{equation}
It follows that either 
\begin{equation*}
\left| f(t)   \right| > |f(st)| - \varepsilon |st|, \qquad \text{or} \qquad\left| f(rt)   \right| > |f(st)| - \varepsilon |st|.
\end{equation*}
 Indeed, otherwise since
\(|a|+|b|\le1\), we would get a contradiction to~\eqref{stIneq}. Since also $\left| \frac t s \right|< \delta$, either 
\begin{equation}\label{trtIneq}
\left| f\left( \frac t s \right)   \right| > |f(t)| - \varepsilon |t|, \qquad \text{or} \qquad \left| f\left( \frac {rt} s \right)   \right| > |f(t)| - \varepsilon |t|
\end{equation}

Now fix \(h\) with \(|h|<\delta\). Repeatedly applying~\ref{trtIneq}
define a sequence \(h_m\) recursively by
\[
h_0=h,
\]
and choosing
\[
h_m= \frac{1}{s}h_{m-1}
\qquad\text{or}\qquad
h_m=\frac{r}{s} h_{m-1}
\]
so that

\begin{equation}\label{ineq1}
|f(h_m)| > |f(h_{m-1})|-\varepsilon |h_{m-1}|.
\end{equation}

In either case, since  \(1<|r|<|s|\),
\begin{equation}\label{ineq2}
0 < |h_m|\le \left|\frac r s  \right||h_{m-1}|.
\end{equation}
where $0 < \left| \frac r s  \right|< 1$.
Therefore, $|h_m|$ strictly decreases to 0, justifying the continued use of \eqref{trtIneq}
Let $c$ be the sum of the geometric series,
\[
c=\sum_{j=0}^\infty \left| \frac r s  \right|^j.
\]
Iterating \eqref{ineq2} gives
\[
0 < |h_j|\le \left| \frac r s  \right|^j|h|
\]
for each non-negative integer $j$ and then iterating \eqref{ineq1} gives
\[
|f(h_m)|
>
|f(h)|-\varepsilon\sum_{j=0}^{m-1}|h_j|
>
|f(h)|-\varepsilon c|h|.
\]
Thus
\[
|f(h)|
<
|f(h_m)|+\varepsilon c|h|.
\]
Since \(h_m\to0\) and \(f(h)=o(1)\), we have \(f(h_m)\to0\). Therefore,
\[
|f(h)|
\le
\varepsilon c|h|.
\]
Because \(\varepsilon>0\) is arbitrary, it follows that
\[
f(h)=o(h).
\]
\end{proof}

\section{The case \texorpdfstring{$|br|+|s|<|a|$}{|br|+|s|<|a|}}\label{sec3}

In this section we consider
$$D(f,h) = af(h)+bf(rh)+f(sh)$$
where $1<|r|<|s|$.  We will show that if $|br|+|s|<|a|$ then $D$ has the MZ-property.

In the rational case, suppose $\frac{\log|r|}{\log|s|}= \frac{k}{j}$ for some
positive integers \(j,k\). Then, this result follows from Corollary~\ref{rat3}. 
In fact, no polynomial $P(z)=\pm z^j \pm bz^k \pm a$
has a root $z$ with $|z|\le |r|^{1/k}$.
To see this, 
 let $\lambda =  |r|^{1/k}$.
We have $\lambda^j = |s|$ and $\lambda^k = |r|$.
If $z$ is a zero of $P(z)$ then $|\pm z^j\pm bz^k|=|a|$. But when
$|z|\le \lambda$,
\[
|\pm z^j\pm bz^k|
\le
\lambda^j+|b|\lambda^k
=
|s|+|b||r|
<
|a|.
\]
So, $P(z)$ has no zeros with $|z|\le \lambda$ and by Corollary~\ref{rat3},
$D(f,h)$ has the MZ-property. 

The following proof works without assuming any rational relation between \(\log |s|\) and
\(\log |r|\).  Furthermore, instead of \(f(h)=o(1)\) it requires only \(f(h)=O(1)\).

\begin{theorem}
Let $1 < |r| < |s|$, and $|br|+|s|<|a|$, and
\[
Df(h)=af(h)+bf(rh)+f(sh).
\]
If $$Df(h) = o(h) \qquad \text{and} \qquad f(h)=O(1)$$
Then
\[
f(h)=o(h).
\]
\end{theorem}
The proof of this theorem is based on repeated applications of the following lemma.
\begin{lemma}
Let
\[
Df(h)=af(h)+bf(rh)+f(sh).
\]
with $1 < |r| < |s|$, and $|br|+|s|<|a|$.
Let $$t=\frac{\log |r|}{\log |s|} \qquad\text{and}\qquad
K = \sum_{k=0}^{\infty}
\left(
\frac{|br|+|s|}{|a|}
\right)^k.
$$
Suppose for some $\alpha \ge 0$
$$Df(h) = o(h) \qquad \text{and} \qquad f(h)=O(h^\alpha).$$
Then for each $\epsilon>0$ there is a $\delta>0$ and  $L>0$ such that whenever $|h|<\delta$,
$$
|f(h)| < K\epsilon |h| + L|h|^{\alpha+t}  
$$
\end{lemma}
We first prove the theorem from the lemma:
\begin{proof}
   Let $1 < |r| < |s|$, and $|br|+|s|<|a|$, and
\[
Df(h)=af(h)+bf(rh)+f(sh).
\]
    Suppose $$Df(h) = o(h) \qquad \text{and} \qquad f(h)=O(1).$$
    We will show that $f(h)=o(h)$.
    
    Since $f(h)=O(h^\alpha)$ for $\alpha=0$.  We apply the lemma with $\epsilon=1$ and get
    positive numbers $\delta_1$ and $L_1$ such that whenever $|h|<\delta_1$, 
    $$
        |f(h)| < K|h| + L_1|h|^{t}.
    $$
    By shrinking $\delta_1$, if necessary, we may assume that $|h|<1$.  Since $t<1$, 
    $|h| < |h|^t$. So, for all $|h|<\delta_1$,
    $$
        |f(h)| < (K+L_1)|h|^{t}.
    $$
    This implies that $f(h)=O(h^t)$.

    We apply the lemma a second time with $\alpha=t$ and $\epsilon=1$.
    We get
    positive numbers $\delta_2$ and $L_2$ such that whenever $|h|<\delta_2$, 
    $$
        |f(h)| < K|h| + L_2|h|^{2t}.
    $$
    If $2t\le 1$ then again, by shrinking $\delta_2$ if necessary, we get
    that for all $|h|<\delta_2$,
    $$
        |f(h)| <  (K + L_2)|h|^{2t},
    $$ 
    and therefore, $f(h)=O(h^{2t})$.

    This process continues as long as $\alpha + t \le 1$.  With each application of the lemma,
    the value of $\alpha$ increases by the constant $t$.  Eventually, after a finite number of applications, we have $f(h)=O(h^{\alpha})$ with $\alpha \le 1$ but $\alpha +t> 1$.  At that point, we apply the lemma one final time, but this time we fix some arbitrary $\epsilon>0$.
    By the lemma, there is a $\delta>0$ and $L>0$ such that whenever $ |h|<\delta$,
    $$
        |f(h)| < K\epsilon |h| + L|h|^{\alpha+t}.
    $$
    Since this time, $\alpha + t > 1$, the first term dominates as $h\to 0$.
    By shrinking $\delta$ if necessary, we can assume $|h|$ is so small that 
    $L|h|^{\alpha+t}<K\epsilon |h| $.  Then we get
    $$
        |f(h)| < 2K\epsilon |h|.
    $$
    Since $\epsilon$ is arbitrary, it follows that $f(h)=o(h)$.
\end{proof}
We now prove the lemma:
\begin{proof}
Let
\[
Df(h)=af(h)+bf(rh)+f(sh).
\]
with $1 < |r| < |s|$, and $|br|+|s|<|a|$.
Let $$t=\frac{\log |r|}{\log |s|} \qquad\text{and}\qquad
K = \sum_{k=0}^{\infty}
\left(
\frac{|br|+|s|}{|a|}
\right)^k.
$$
Fix $\alpha \ge 0$ such that 
$$Df(h) = o(h) \qquad \text{and} \qquad f(h)=O(h^\alpha).$$
Fix $\epsilon>0$.  
We must find $\delta>0$ and  $L>0$ such that whenever $|h|<\delta$,
$$
|f(h)| < K\epsilon |h| + L|h|^{\alpha+t}.  
$$
Using that $f(h)=O(h^\alpha)$, choose $\delta>0$ and $B>0$, such that 
\begin{equation}\label{defB}
|f(h)|\le B|h|^\alpha
\qquad\text{for all } |h|<|s|\delta .
\end{equation}
Since
\[
af(h)+bf(rh)+f(sh)=o(h),
\]
shrink $\delta$, if necessary, so that for all $|h|<\delta$,
\begin{equation}\label{est1}
|f(h)|
\le
\frac1{|a|}|f(sh)|
+
\frac{|b|}{|a|}|f(rh)|
+
\epsilon |h|.
\end{equation}
If $|sh|<\delta$ (and therefore $|rh|<\delta$), then applying the same estimate to $rh$ and $sh$ gives
\begin{equation*}
|f(rh)|
\le
\frac1{|a|}|f(rsh)|
+
\frac{|b|}{|a|}|f(r^2h)|
+
\epsilon |rh|,
\end{equation*}
and
\begin{equation*}
|f(sh)|
\le
\frac1{|a|}|f(s^2h)|
+
\frac{|b|}{|a|}|f(rsh)|
+
\epsilon |sh|.
\end{equation*}
Substituting these into~\eqref{est1}, we obtain
\begin{equation}\label{est2}
|f(h)|
\le
\frac1{|a|^2}
\left(
|f(s^2h)|
+
2|b||f(rsh)|
+
|b|^2|f(r^2h)|
\right)
+
\epsilon |h|
\left(
1+\frac{|br|+|s|}{|a|}
\right).
\end{equation}
If also $|s^2h|<\delta$, then applying ~\eqref{est1} to $|r^2h|,|rsh|,|s^2h| $ we obtain
\begin{eqnarray*}
|f(r^2h)| &\le& \frac1{|a|}|f(r^2sh)|+\frac{|b|}{|a|}|f(r^3h)|+\epsilon |r^2h|,\\
|f(rsh)| &\le& \frac1{|a|}|f(rs^2h)|+\frac{|b|}{|a|}|f(r^2sh)|+\epsilon |rsh|,\\
|f(s^2h)| &\le& \frac1{|a|}|f(s^3h)|+\frac{|b|}{|a|}|f(rs^2h)|+\epsilon |s^2h|.
\end{eqnarray*}
Substituting into~\eqref{est2}, we obtain
\begin{eqnarray*}
|f(h)|&\le&\frac1{|a|^3}\left(|f(s^3h)|+3|bf(rs^2h)|+3|b^2f(r^2sh)|+|b^3f(r^3h)|\right)\\
&&+ \epsilon |h|\left(
1+\frac{|br|+|s|}{|a|}+\left(\frac{|br|+|s|}{|a|}\right)^2
\right).
\end{eqnarray*}
Continuing inductively, if $|s^Nh|<\delta,$
then
\begin{equation}\label{estN}
|f(h)|
\le
\frac1{|a|^{N+1}}
\sum_{k=0}^{N+1}
\binom{N+1}{k}
|b^k f(r^ks^{N+1-k}h)|
+
\epsilon |h|
\sum_{k=0}^{N}
\left(
\frac{|br|+|s|}{|a|}
\right)^k .   
\end{equation}
Furthermore, since
\[
|r^ks^{N+1-k}h|
\le
|s^{N+1}h|
<
|s|\delta,
\]
it follows from \eqref{defB} that 
\[
|f(r^ks^{N+1-k}h)|\le B|h|^\alpha.
\]
Hence
\begin{equation*}|f(h)|\le\frac {B|h|^\alpha} {|a|^{N+1}}\sum_{k=0}^{N+1}\binom{N+1}{k}|b|^k
+
\epsilon |h|\sum_{k=0}^{N}\left(\frac{|br|+|s|}{|a|}\right)^k    
\end{equation*}
Applying the definition of $K$ and the binomial theorem,
\[
\sum_{k=0}^{N+1}
\binom{N+1}{k}|b|^k
=
(1+|b|)^{N+1},
\]
we get
\begin{equation}\label{mainEq}
|f(h)|<
B|h|^\alpha\left(
\frac{1+|b|}{|a|}
\right)^{N+1}
+
K\epsilon |h|.
\end{equation}
Now fix $N$ such that
\[
|s^Nh|<\delta\le |s^{N+1}h|.
\]
Then,
\begin{eqnarray*}
    \left(\frac{1+|b|}{|a|}\right)^{N+1} &<& \left(\frac{|s|+|br|}{|ar|}\right)^{N+1}  \\
&<& \left(\frac{|a|}{|ar|}\right)^{N+1} \\
&=& \frac1{|r|^{N+1}}\\
&=& \frac1{|s|^{t(N+1)}}.
\end{eqnarray*}
Since
\[
\delta\le |s|^{N+1}|h|,
\]
it follows that
\[
\frac1{|s|^{t(N+1)}}
\le
\left(
\frac{|h|}{\delta}
\right)^{t},
\]
yielding 
$$
 \left(\frac{1+|b|}{|a|}\right)^{N+1} < \left(\frac{|h|}{\delta} \right)^{t}.
$$
Substituting into \eqref{mainEq} and letting $L =\frac B {\delta^t} $, we obtain
\begin{eqnarray*}
|f(h)|
&\le&
B|h|^\alpha \left( \frac {|h|} \delta\right)^{t}
+
K\epsilon |h|\\
&=& L|h|^{\alpha +t}
+
K\epsilon |h|
\end{eqnarray*}
\end{proof}

\section{The case \texorpdfstring{$|a|+1\le |b|$ and $|s|+|a|<|br|$}{|a|+1 <= |b| and $|s|+|a|< |br|$}}\label{sec4}

In this section we prove the following theorem.

\begin{theorem}\label{wholeTheorem}
Suppose that
$$
1<|r|<|s|\qquad\text{and}\qquad|a|+1\leq |b|\qquad\text{and}\qquad|s| + |a|< |br|.
$$
Then the difference
$$
D(f,h) = af(h)+bf(rh)+f(sh)
$$
has the MZ-property.
\end{theorem}

This time, we do not have a proof that handles both the rational and irrational cases simultaneously. We will therefore split them.  First, the rational case:

\begin{theorem}\label{ratCase}
Suppose that
$$
1<|r|<|s|\qquad\text{and}\qquad|a|+1\leq |b|\qquad\text{and}\qquad|s| + |a|< |br|.
$$
Suppose also that $\frac{\log |r|}{\log |s|}\in\mathbb Q$,
Then the difference
$$
D(f,h) = af(h)+bf(rh)+f(sh)
$$
has the MZ-property.
\end{theorem}

\begin{proof}
Let $\frac{\log |r|}{\log |s|} = \frac k j$ where $k,j$ are positive integers.
Let $\lambda = |s|^{1/j}$, so that $\lambda^j = |s|$ and $\lambda^k = |r|$.  By Corollary~\ref{rat3}, it is enough to show that no trinomial of the form
$$
P(z) = \pm z^j \pm  bz^k \pm a
$$
can have zeros in the annulus $1<|z|\leq \lambda$.

If $P(z)=0$ then $|bz^k|=|\pm z^j\pm a| \le |z|^j + |a|$. In other words, $\phi(x) = x^j - |b|x^k + |a| \ge 0$ when $x=|z|$.  We will show this can't happen for $x\in (1,\lambda]$.
By Descartes Rule of Signs,
$\phi$ has at most two positive real roots.  To estimate them,
$$
\phi(0)=|a|>0.
$$
\[
\phi(1)=1-|b|+|a|\le0
\]
\[
\phi(\lambda)=|s|-|b||r|+|a|<0.
\]
and for large enough $x$,
\[
\phi(x) > 0.
\]
So, $\phi(x)$ has a root in $(0,1]$ and another root in $(\lambda,\infty)$ and is strictly
negative in the region $(1,\lambda]$. 
\end{proof}

We now consider the same theorem for the irrational case.
To prove this theorem, it will be convenient to use the following lemma.
Examples~\ref{Ex1} and \ref{Ex2} demonstrate
that this next lemma is false in the rational case, explaining why we split
into cases.

\begin{lemma}\label{falseLemma}
Let $$D_1(f,h)=af(h)+bf(rh)+f(sh)$$ and $$D_2(f,h)=cf(h)+df(rh)+f(sh)$$ be two irrational differences with the same nodes.
Suppose $|a|=|c|$ and $|b|=|d|$.  Then $D_1$ has the MZ-property if and only if $D_2$ has the MZ-property.
\end{lemma}

\begin{proof}
Suppose that $D_1$ does not have the MZ-property. Then for some function $f:\mathbb R \to \mathbb C$ we have that $D_1(f,h)=o(h)$, $f=o(1)$, but $f\neq o(h)$.  We will finish the proof by finding a function $g:\mathbb R \to \mathbb C$ such that $D_2(g,h)=o(h)$, $g=o(1)$, but $g\neq o(h)$.  

Since $|a|=|c|$, choose \(\theta_1\in\mathbb R\) such that
\(
c=ae^{i\theta_1}.
\)
Since $|b|=|d|$, choose \(\theta_2\in\mathbb R\) such that
\(
de^{i\theta_2}=be^{i\theta_1}.
\)

Define
\[
H=\{r^ns^m:n,m\in\mathbb Z\}.
\]
Since
\(
\frac{\log |r|}{\log |s|}\notin\Q,
\)
every element of \(H\) has a unique representation of the form
\(
r^ns^m.
\)

Using that $f$ is not $o(h)$, 
let \(h_k\to0\) be such that \(f(h_k)\ne o(h_k)\). 
Consider the collection of cosets
$
\{h_kH\},
$
and note that this collection may have duplicates.
If $h\in \bigcup h_kH$
we will define $g(h)$ as follows.
Find the first $k$ such that $h\in h_kH$, and express $h$ as $h=h_kr^ns^m$.
Then let $g(h)=e^{i(m\theta_1+n\theta_2)}f(h)$.
If $h\notin \bigcup h_kH$, simply let $g(h)=0$.

Since $|g(h)|\le|f(h)|$, we get $g=o(1)$.
Since $|g(h)|=|f(h)|$, for $h\in \bigcup h_kH$, we get  
$g(h_k)\ne o(h_k)$.  Therefore, $g\neq o(h)$.
It remains to show that $D_2(g,h)=o(h)$.

If $h\in \bigcup h_kH$,
then the first $k$ such that $h\in h_kH$ is also the first $k$ such that $rh\in h_kH$
and the first $k$ such that $sh\in h_kH$.  Then writing $h=h_kr^ns^m$ we have 
\[
g(h) = e^{i(m\theta_1+n\theta_2)}f(h),
\qquad
g(rh)=e^{i\theta_2}e^{i(m\theta_1+n\theta_2)}f(rh),
\qquad
g(sh)=e^{i\theta_1}e^{i(m\theta_1+n\theta_2)}f(sh).
\]
Hence,
\[
cg(h)+dg(rh)+g(sh)=e^{i(m\theta_1+n\theta_2)}\left(cf(h)+de^{i\theta_2}f(rh)+e^{i\theta_1}f(sh)\right).
\]
By choice of \(\theta_1,\theta_2\), 
\[
cf(h)+de^{i\theta_2}f(rh)+e^{i\theta_1}f(sh)=e^{i\theta_1}\left(af(h)+bf(rh)+f(sh)\right),
\]
and therefore
\begin{eqnarray*}
    D_2(g,h) &=& cg(h)+dg(rh)+g(sh)\\
    &=& e^{i\theta_1}\left(af(h)+bf(rh)+f(sh)\right)\\
    &=& e^{i\theta_1}D_1(f,h) = o(h).
\end{eqnarray*}
\end{proof}

This lemma will be combined with the following theorem.
Although it is only needed in the irrational case, the proof of the next theorem goes through
for either rational or irrational cases.

\begin{theorem}\label{tgr}
Let $c, d$ be real numbers such that
\[
0<|c|<1<|d|,
\]
and let $c_1>0, d_1>0$ be real numbers such that 
\[
c_1+d_1\le1 \qquad\text{and} \qquad c_1|c|+d_1|d|< 1.
\]
Then
\[
D(f,h) = f(h)-c_1f(ch)-d_1f(dh)
\]
has the MZ-property.
\end{theorem}
\begin{proof}
Set
\[
\lambda:=\frac{1-c_1|c|}{|d|d_1} \qquad\text{and}\qquad \rho:=\frac{\log \lambda}{\log |d|}.
\]
By hypothesis, \(\lambda>1\) and $\rho>0$.

Suppose $D(f,h)=o(h)$ and $f(h)=o(1)$.
We will show $f(h)=o(h)$.
Suppose on the contrary, 
\[
\frac{f(h)}{h}\not\to 0.
\]
Replacing \(f\) with \(-f\) if necessary, we may assume that there exists \(\alpha>0\) such that
for arbitrarily small \(h\),
\[
\frac{f(h)}{|h|}>\alpha.
\]
Since \(f(h)=o(1)\) and $D(f,h)=o(h)$,
there exists \(\beta>0\) such that whenever \(|h|<\beta\) we have
\[
|f(h)|<1
\]
and
\[
|D(f,h)| = |f(h)-c_1f(ch)-d_1f(dh)|
<
\frac{(\lambda-1)|d|d_1}{4}\alpha |h|.
\]

Choose \(h_0\) small enough that the following three properties hold:
\[
|h_0|<\frac{\beta}{|d|},
\qquad
\frac{f(h_0)}{|h_0|}>\alpha,
\qquad\text{and}\qquad
|h_0|^\rho<\left(\frac{\beta}{|d|}\right)^{1+\rho}\frac{\alpha}{4}.
\]
The last inequality is possible because \(\rho>0\), so \(|h|^\rho\to 0\) as \(h\to 0\).
Now, set
\[
\epsilon_0=\frac{f(h_0)}{2} > 0,
\]
and define
\[
G(h):=\frac{f(h)-\epsilon_0}{|h|}.
\]
We define a sequence \(h_0,h_1,h_2,\dots\) recursively as follows. For each \(i\ge 0\):

\medskip

\noindent
\textbf{Case 1.} If
\[
f(ch_i)\ge |c|f(h_i)+(1-|c|)\epsilon_0,
\]
then
\[
h_{i+1}=ch_i.
\]

\noindent
\textbf{Case 2.} If
\[
f(ch_i)<|c|f(h_i)+(1-|c|)\epsilon_0,
\]
then
\[
h_{i+1}=dh_i.
\]

\medskip

We now state and prove several claims.
\medskip

\noindent
\textbf{Claim 0.} 
$G(h_0) >\frac{\alpha}{2}$.\newline
Indeed,
\[
G(h_0)=\frac{f(h_0)-\epsilon_0}{|h_0|}
=\frac{f(h_0)/2}{|h_0|}
>\frac{\alpha}{2}.
\]

\bigskip
\noindent
\textbf{Claim 1.} In Case 1, $G(h_{i+1})\ge G(h_i)$.\newline
Indeed, 
\[
G(h_{i+1})
=
\frac{f(ch_i)-\epsilon_0}{|c||h_i|}
\ge
\frac{|c|f(h_i)+(1-|c|)\epsilon_0-\epsilon_0}{|c||h_i|}
=
\frac{f(h_i)-\epsilon_0}{|h_i|}
=
G(h_i).
\]

\bigskip
\noindent
\textbf{Claim 2.} In Case 2, if \(|h_i|<\beta\). Then
\[
G(h_{i+1})
>
\lambda G(h_i)-\frac{\lambda-1}{4}\alpha.
\]

\smallskip
\noindent
\emph{Proof of Claim 2:}
Since \(|h_i|<\beta\),
\[
|f(h_i)-c_1f(ch_i)-d_1f(dh_i)|
<
\frac{(\lambda-1)|d|d_1}{4}\alpha |h_i|.
\]
Therefore,
\[
d_1f(dh_i)
>
f(h_i)-c_1f(ch_i)-\frac{(\lambda-1)|d|d_1}{4}\alpha |h_i|.
\]
Subtracting \(d_1\epsilon_0\) from both sides and dividing by \(|d||h_i|\), we get
\[
d_1G(h_{i+1})
=
d_1\frac{f(dh_i)-\epsilon_0}{|d||h_i|}
\]
\[
>
\frac{
f(h_i)-c_1f(ch_i)-\frac{(\lambda-1)|d|d_1}{4}\alpha |h_i|-d_1\epsilon_0
}{|d||h_i|}.
\]
Now Case 2 means that
\[
f(ch_i)<|c|f(h_i)+(1-|c|)\epsilon_0,
\]
so
\[
-c_1f(ch_i)>-c_1|c|f(h_i)-c_1(1-|c|)\epsilon_0.
\]
Substituting this gives
\[
d_1G(h_{i+1})
>
\frac{
(1-c_1|c|)f(h_i)-\bigl(c_1(1-|c|)+d_1\bigr)\epsilon_0-\frac{(\lambda-1)|d|d_1}{4}\alpha |h_i|
}{|d||h_i|}.
\]
Since \(c_1+d_1\le 1\), we have
\[
c_1(1-|c|)+d_1 \le 1-c_1|c|.
\]
Therefore,
\[
d_1G(h_{i+1})
>
\frac{
(1-c_1|c|)\bigl(f(h_i)-\epsilon_0\bigr)-\frac{(\lambda-1)|d|d_1}{4}\alpha |h_i|
}{|d||h_i|}.
\]
Hence,
\[
d_1G(h_{i+1})
>
\frac{1-c_1|c|}{|d|}G(h_i)-\frac{(\lambda-1)d_1}{4}\alpha,
\]
and dividing by \(d_1\) yields the claim.

\bigskip
\noindent
\textbf{Claim 3.} In Case 2, if \(|h_i|<\beta\) and  \(G(h_i)\ge \alpha/2\) then
\(G(h_{i+1})> \alpha/2\).\newline
Indeed, if \(G(h_i)\ge \alpha/2\), then by Claim 2,
\[
G(h_{i+1})
>
\lambda\frac{\alpha}{2}-\frac{\lambda-1}{4}\alpha
=
\frac{\lambda+1}{4}\alpha
>
\frac{\alpha}{2},
\]
since \(\lambda>1\).

\bigskip
\noindent
\textbf{Claim 4.} If $\max\{h_0, h_1 , \ldots, h_n\}<\beta$
then $\min\{f(h_0), f(h_1), \ldots, f(h_n)\}>\epsilon_0$.

\smallskip
\noindent
\emph{Proof of Claim 4:}  
By Claim 0, 
\[
G(h_0)>\frac{\alpha}{2}.
\]
Suppose for some $i\le n$,
\[
G(h_i)>\frac{\alpha}{2}.
\]
Then, using Claim 1 in Case 1, or Claim 3 in Case 2, 
\[
G(h_{i+1})>\frac{\alpha}{2}.
\]
Therefore, by induction,
\[
G(h_i)>\frac{\alpha}{2}
\qquad\text{for all \(i = 0, 1, \ldots, n\) }.
\]
Consequently,
\[
f(h_i)-\epsilon_0=|h_i|G(h_i)>0,
\]
so
\[
f(h_i)>\epsilon_0
\qquad\text{for all \(i = 0, 1, \ldots, n\) }.
\]

\bigskip
\noindent
\textbf{Claim 5.} 
If $\max\{h_0, h_1 , \ldots, h_n\}<\beta$
and Case 2 is used $m$ times in the sequence $\{h_0, h_1 , \ldots, h_n\}$
then $$G(h_n)> \frac{\lambda^m\alpha}{4}.$$\newline

\smallskip
\noindent
\emph{Proof of Claim 5:}  
If \(m\) be the number of times Case 2 is used in the construction of
\[
h_0,h_1,\dots,h_n.
\]
Then Case 1 is used exactly \(n-m\) times, so
\[
|h_n|=|h_0||c|^{\,n-m}|d|^m.
\]
Each time Case 1 occurs, Claim 1 shows that \(G\) does not decrease. 
\[
T(x):=\lambda x-\frac{\lambda-1}{4}\alpha.
\]
Since Case 2 occurs exactly \(m\) times, we obtain from Claim 2 that 
\[
G(h_n)> T^m(G(h_0)).
\]
Expanding this recursion gives
\[
G(h_n)> \lambda^mG(h_0)-\frac{\lambda-1}{4}\alpha\left(1+\lambda+\cdots+\lambda^{m-1}\right).
\]
Since
\[
1+\lambda+\cdots+\lambda^{m-1}=\frac{\lambda^m-1}{\lambda-1},
\]
this becomes
\[
G(h_n)> \lambda^mG(h_0)-\frac{\lambda^m-1}{4}\alpha
>
\lambda^m\left(G(h_0)-\frac{\alpha}{4}\right).
\]
Since \(G(h_0)>\alpha/2\) by Claim 0, we conclude that
\[
G(h_n)> \lambda^m\frac{\alpha}{4}.
\]

\bigskip
\noindent
\textbf{Claim 6.} 
For all \(i\),
\[
|h_i| <\beta.
\]

\smallskip
\noindent
\emph{Proof of Claim 6:}  
In fact, we will prove the stronger inequality $|h_i| < \frac{\beta}{|d|}$.
Suppose not. Then there exists a smallest \(i\ge 0\) such that
\[
|h_i|\ge\frac{\beta}{|d|}.
\]
Since \(|h_0|<\frac{\beta}{|d|}\) by construction, we must have \(i\ge 1\).

By minimality of \(i\), we have
\[
|h_{i-1}| < \frac{\beta}{|d|}.
\]
Since \(h_i\) is obtained from \(h_{i-1}\) by multiplying either by \(c\) or by \(d\), we have
\[
|h_i|=|c||h_{i-1}| \quad \text{or} \quad |h_i|=|d||h_{i-1}|.
\]
In either case, using \(|h_{i-1}| < \frac{\beta}{|d|}\) and \(0<|c|<1<|d|\), we obtain
\[
|h_i| < \beta. 
\]
Let $m$ be the number of times we used Case 2 in the sequence $\{h_0,h_1, \ldots, h_i\}$.  Since
\[
|h_i|=|h_0||c|^{\,i-m}|d|^m\le |h_0||d|^m,
\]
we have
\[
|d|^m\ge \frac{|h_i|}{|h_0|}.
\]
Because
\[
\rho=\frac{\log\lambda}{\log|d|},
\]
we have
\[
\lambda=|d|^\rho,
\]
and therefore
\[
\lambda^m=|d|^{\rho m}\ge \left(\frac{|h_i|}{|h_0|}\right)^\rho.
\]
Using Claim 5, we obtain
\[
f(h_i)-\epsilon_0
=
|h_i|G(h_i)
>
|h_i|\lambda^m\frac{\alpha}{4}
\ge
|h_i|\left(\frac{|h_i|}{|h_0|}\right)^\rho\frac{\alpha}{4}.
\]
Since \(|h_i|\ge\beta/|d|\), this yields
\[
f(h_i)-\epsilon_0
\ge
\frac{\alpha}{4}\left(\frac{\beta}{|d|}\right)^{1+\rho}|h_0|^{-\rho}.
\]
By the choice of \(h_0\),
\[
|h_0|^\rho<\frac{\alpha}{4}\left(\frac{\beta}{|d|}\right)^{1+\rho},
\]
and therefore
\[
\frac{\alpha}{4}\left(\frac{\beta}{|d|}\right)^{1+\rho}|h_0|^{-\rho}>1.
\]
Hence,
\[
f(h_i)-\epsilon_0>1.
\]
Since \(\epsilon_0>0\), we obtain
\[
f(h_i)>1.
\]
But this contradicts the choice of \(\beta\), because  \(|h_i|<\beta\) and because
\(|f(h)|<1\) whenever \(|h|<\beta\).
This contradiction establishes Claim 6.

\bigskip
\noindent
\textbf{Claim 7.} $\epsilon_0 < f(h_i) < 1$ for all \(i\).

\smallskip
\noindent
\emph{Proof of Claim 7:} 
By Claim 6, \(|h_i| < \beta\) for all \(i\).  Therefore, by choice of $\beta$,
$f(h_i)<1$ for all $i$, and 
by Claim 4, $f(h_i)>\epsilon_0$ for all \(i\).

\bigskip
\noindent
\textbf{Claim 8.} \(|h_i|\to 0\).

\smallskip
\noindent
\emph{Proof of Claim 8:}
Suppose first that Case 2 occurs infinitely many times. Then Claims 1 and 5 imply 
\(G(h_i)\to\infty\).  By definition, $G(h_i)=\frac{f(h_i)-\epsilon_0}{|h_i|}$ and by Claim 7, \(|f(h_i)|<1\).  So, this is possible only if \(|h_i|\to 0\). 

Now suppose that Case 2 occurs only finitely many times. Then from some point onward only Case 1 occurs, so
\[
h_{i+1}=ch_i
\]
for all sufficiently large \(i\). Hence
\[
|h_{i+1}|=|c|\,|h_i|
\]
for all sufficiently large \(i\), and since \(0<|c|<1\), it follows that \(h_i\to 0\).
This finishes the proof of Claim 8.

\bigskip
\noindent
To finish the proof of the theorem, note that Claims 7 and 8 directly contradict that \(f(h)=o(1)\).
Therefore our assumption was false, and hence
\[
f(h)=o(h).
\]
\end{proof}

We are now ready to handle the irrational case.  This next theorem together with Theorem~\ref{ratCase} establishes Theorem~\ref{wholeTheorem}.

\begin{theorem}
Let $D(f,h) = af(h)+bf(rh)+f(sh)$ be an irrational difference with $1<|r|<|s|$.  
Suppose that
$$|a|+1\leq |b|\qquad\text{and}\qquad|s| + |a|< |br|.
$$
Then $D(f,h)$ has the MZ-property.
\end{theorem}

\begin{proof}
Let 
$$
c = \frac 1 r, 
\qquad  d = \frac s r, 
\qquad c_1 = \left|\frac a b\right|, 
\qquad \text{ and } d_1 = \left|\frac 1 b\right|.  
$$
Then 
\begin{eqnarray*}
1<|r|<|s| &\text{implies}& 0<|c|<1<|d|,\\
|a|+1\leq |b| &\text{implies}&  c_1+d_1\le1,\\
|s| + |a|< |br| &\text{implies}&  \left|\frac s {br}\right| + \left|\frac a {br}\right|< 1,\qquad\text{ i.e. } \qquad c_1|c|+d_1|d|< 1  
\end{eqnarray*}
Therefore, by Theorem~\ref{tgr}, 
\begin{eqnarray*}
f(h)-c_1f(ch)-d_1f(dh)=-c_1f(ch)+f(h)-d_1f(dh)
\end{eqnarray*}
has the MZ-property.  Recall that the MZ-property is invariant under scaling the nodes and coefficients.
Dividing the nodes by $c$ shows that 
$$-c_1f(h)+f(\frac 1 c h)-d_1f(\frac d c h)
= -c_1f(h)+f(r h)-d_1f(s h)
$$ has the MZ-property.  Multiplying the coefficients by $|b|$ shows that 
$$ -|b|c_1f(h)+|b|f(r h)-|b|d_1f(s h) = -|a|f(h)+ |b| f(r h)- f(s h)$$ has the MZ-property.
Then, applying Lemma~\ref{falseLemma},  we get that
$$D(f,h) = af(h)+bf(rh)+f(sh)$$ has the MZ-property.
\end{proof}


\section{The cases where MZ-property fails}\label{sec5}
\begin{theorem}\label{tpqf}
Let $r,s$ be nonzero real numbers with $\frac{\log  |r|}{\log |s|}\notin\mathbb{Q}$,
and let $a,b\in\mathbb{C}$. Assume that
\[
|a|+|b|>1
\qquad\text{and}\qquad
\bigl||a|-|br|\bigr|\le |s|.
\]
Then for the difference
\[
D(f,h) = af(h)+bf(rh)+f(sh)
\]
there exists a complex-valued function $f$ such that

\[
f(h)=o(1),\quad f(h)\notin o(h),\quad\text{and}\quad D(f,h)\equiv 0
\text{ for all }h\in\mathbb{R}.
\]
\end{theorem}

\begin{proof}
Since
\[
\max_{\theta\in\mathbb R}\bigl|-a+|b|e^{i\theta}\bigr|=|a|+|b|>1,
\]
there exists $\tau\in\mathbb R$ such that
\[
\bigl|-a+|b|e^{i\tau}\bigr|>1.
\]
Also,
\[
\min_{\theta\in\mathbb R}\bigl|-a+|b||r|e^{i\theta}\bigr|
=\bigl||a|-|b||r|\bigr|
\le |s|,
\]
so there exists $\sigma\in\mathbb R$ such that
\[
\bigl|-a+|b||r|e^{i\sigma}\bigr|\le |s|.
\]

Define
\[
\Psi(t,\theta)=|s|^t-\bigl|-a+|b||r|^te^{i\theta}\bigr|,
\qquad (t,\theta)\in[0,1]\times\mathbb R.
\]
Then
\[
\Psi(0,\tau)<0
\qquad\text{and}\qquad
\Psi(1,\sigma)\ge 0.
\]

Let
\[
\theta(t)=(1-t)\tau+t\sigma
\qquad\text{and}\qquad
F(t)=\Psi(t,\theta(t)).
\]
Then \(F\) is continuous on \([0,1]\), with
\[
F(0)<0
\qquad\text{and}\qquad
F(1)\ge 0.
\]
Hence, by the intermediate value theorem, there exists \(t_0\in(0,1]\) such that
\[
F(t_0)=0.
\]
Thus
\[
|s|^{t_0}=\bigl|-a+|b||r|^{t_0}e^{i\theta(t_0)}\bigr|.
\]
Therefore there exists \(\theta_1\in\mathbb R\) such that
\[
|s|^{t_0}e^{i\theta_1}
=
-a+|b||r|^{t_0}e^{i\theta(t_0)},
\]
and \(\theta_2\in\mathbb R\) such that 
\[
a + b\,|r|^{t_0}e^{i\theta_2}+ |s|^{t_0}e^{i\theta_1}=0.
\tag{1}\label{phase-relation}
\]
Let
\[
H=\{|r|^n|s|^m:\ n,m\in\mathbb Z\}.
\]
Since \(\log(|r|)/\log(|s|)\notin\mathbb Q\), each \(h\in H\) has a unique representation
\[
h=|r|^n|s|^m.
\]
Define
\[
f(h)=
\begin{cases}
e^{i(n\theta_2+m\theta_1)}|h|^{t_0}, & h\in H,\ h=|r|^n|s|^m,\\[4pt]
0, & h\notin H.
\end{cases}
\]
Since $0<t_0\le 1$, $f(h) = o(1)$ and $f(h)\ne o(h)$.
For \(h\in H\), and hence both \(rh\in H\) and \(sh\in H\),
\[
D(f,h)
=
\bigl(a+b\,|r|^{t_0}e^{i\theta_2}+|s|^{t_0}e^{i\theta_1}\bigr)f(h)=0
\]
by \eqref{phase-relation}.
For \(h\notin H\), and hence both \(rh\notin H\) and \(sh\notin H\), 
\[
D(f,h)=0.
\]
Thus
\[
D(f,h)\equiv 0
\qquad\text{for all }h\in\mathbb R.
\]
\end{proof}

The following examples show that Theorem~\ref{tpqf} can go either way in the rational case.

\begin{example}\label{Ex4}
Let $D(f,h) = f(h)+2f(4h)+f(8h)$, with $a=1, b=2, r=4, s=8$.  Then $D(f,h)$ is rational with $|a|+|b|>1$ and
$||a|-|br||=7<|s|$. The inequalities of Theorem~\ref{tpqf} are satisfied. Since the nodes are positive, there is only one characteristic polynomial, $x^3+2x^2+1$.  
The roots are approximately $-2.21$, $0.1+0.67i$, and $0.1-0.67i$.
The absolute values of its roots are 2.21 and 0.67, which are outside of
$(1,2]$ so by Theorem~\ref{rat3}, $D(f,h)$ has the MZ-property.
\end{example}

\begin{example}\label{Ex5}
Let $D(f,h) = f(h)-3f(2h)+f(8h)$.  Then $D(f,h)$ is rational with $a=1, b=-3, r=2, s=8$.  Then $|a|+|b|>1$ and
$||a|-|br||=5<|s|$. The inequalities of Theorem~\ref{tpqf} are satisfied. Since the nodes are positive, there is only one characteristic polynomial, $x^3-3x+1$.  One of its roots, 1.53, is in the region
$(1,2]$. So, by Theorem~\ref{rat3}, $D(f,h)$ does not have the MZ-property.
\end{example}

\begin{theorem}\label{sbra}
Suppose that
\[
1+|a|>|b|
\qquad\text{and}\qquad
|s|<|b||r|-|a|.
\]
Then for the difference
\[
D(f,h)=|a|f(h) -|b|f(rh)+f(sh)
\]
there exists a complex-valued function $f$ such that
\[
f(h)=o(1),\quad f(h)\notin o(h),\quad\text{and}\quad D(f,h)\equiv 0
\text{ for all }h\in\mathbb{R}.
\]
\end{theorem}

\begin{proof}
Define
\[
g(t)=|s|^{\log_{|r|}t}-|b|t+|a|.
\]
Then
\[
g(1)=1-|b|+|a|>0
\]
and
\[
g(|r|)=|s|-|b||r|+|a|<0.
\]
Hence, by continuity, there exists
\[
1<t_0<|r|
\]
such that
\[
g(t_0)=0.
\]

Define
\[
f(h)=|h|^{\log_{|r|}t_0}.
\]
Then
\[
f(sh)=|s|^{\log_{|r|}t_0}f(h)
\]
and
\[
f(rh)=t_0f(h).
\]
Therefore
\[
D(f,h)
=
\left(
|a|
-|b|t_0
+
|s|^{\log_{|r|}t_0}
\right)f(h)
=
g(t_0)f(h)
=
0.
\]
Since
\[
1<t_0<|r|,
\]
we have
\[
0<\log_{|r|}t_0<1,
\]
and consequently
\[
f(h)= o(1)\qquad\text{and}\qquad f(h)\neq o(h).
\]
\end{proof}

We now complete the proof of Theorem~\ref{mainTheorem}.
\begin{theorem}
Let $D(f,h) = af(h)+bf(rh)+f(sh)$ be an irrational difference with 
    distinct nodes, $1,r,s$ where $1\le |r|\le |s|$.  If $D$ has the MZ-property then it 
    satisfies at least one of the following conditions:
    \begin{itemize}
\item $|a|+|b|\leq 1$
\item $|br| + |s| < |a|$
\item $|a|+1\leq |b|$ and $|s| + |a|< |br|$.
    \end{itemize}
\end{theorem}

\begin{proof}
Suppose $D(f,h) = af(h)+bf(rh)+f(sh)$ is irrational and does not satisfy any of the three conditions. We will show it does not have the MZ-property. Since $D(f,h)$ is irrational, we have
$1<|r|<|s|$.

Case 1: $|a|+|b|>1$, $|br|+|s|\ge |a|$, and $|s|+|a|\ge |br|$.
Then $|s|\ge |br|-|a|$ and $|s|\ge |a|-|br|$. So, $|s|\ge ||a|-|br||$.  By Theorem~\ref{tpqf}, $D(f,h)$ does not have the MZ-property.

Case 2: $|a|+|b|>1$, $|br|+|s|\ge |a|$, and $|s|+|a|< |br|$.
Then $|a|+1>|b|$, so by Theorem~\ref{sbra}, the difference  $|a|f(h) -|b|f(rh)+f(sh)$ fails the MZ-property.
Then by Lemma~\ref{falseLemma}, $D(f,h) = af(h)+bf(rh)+f(sh)$ also fails the MZ-property.

\end{proof}

\end{document}